\documentclass[preprint,12pt]{elsarticle}

\usepackage[letterpaper, margin = 2cm]{geometry} 

\usepackage{amssymb}
\usepackage{amsmath}
\usepackage{lipsum}
\usepackage{amsfonts}
\usepackage{graphicx,subfigure}
\usepackage{epstopdf}
\usepackage{algorithmic}
\usepackage[ruled,linesnumbered]{algorithm2e}
\usepackage{booktabs}
\usepackage{multirow} 
\usepackage{enumitem}
\usepackage{stmaryrd}
\usepackage{bm}

\usepackage{hyperref}
\usepackage{xcolor}
\hypersetup{
    colorlinks=true,
    citecolor=green, 
    linkcolor=red,   
    urlcolor=cyan    
}

\newcommand{\cu}{\boldsymbol}

\newproof{lemma}{Lemma}
\newtheorem{remark}{Remark}

\makeatletter
\def\ps@pprintTitle{%
 \let\@oddhead\@empty
 \let\@evenhead\@empty
 \let\@oddfoot\@empty
 \let\@evenfoot\@empty
}
\makeatother

\begin{document}

\begin{frontmatter}

\title{Momentum-Accelerated Richardson(m) and Their Multilevel Neural Solvers}

\author[a]{Zhen Wang}
\author[a]{Yun Liu}
\author[a]{Chen Cui}
\ead{cuichensx@gmail.com}
\author[a]{Shi Shu}
\ead{shushi@xtu.edu.cn}

\affiliation[a]{organization={Hunan Key Laboratory for Computation and Simulation in Science and Engineering,
Key Laboratory of Intelligent Computing and Information Processing of Ministry
of Education, School of Mathematics and Computational Science, Xiangtan University},
            city={Xiangtan},
            postcode={411105}, 
            state={Hunan},
            country={China}}

\begin{abstract}
    Recently, designing neural solvers for large-scale linear systems of equations has emerged as a promising approach in scientific and engineering computing. This paper first introduce the Richardson(m) neural solver by employing a meta network to predict the weights of the long-step Richardson iterative method. Next, by incorporating momentum and preconditioning techniques, we further enhance convergence.
    Numerical experiments on anisotropic second-order elliptic equations demonstrate that these new solvers achieve faster convergence and lower computational complexity compared to both the Chebyshev iterative method with optimal weights and the Chebyshev semi-iteration method.
    To address the strong dependence of the aforementioned single-level neural solvers on PDE parameters and grid size, we integrate them with two multilevel neural solvers developed in recent years. Using alternating optimization techniques, we construct Richardson(m)-FNS for anisotropic equations and NAG-Richardson(m)-WANS for the Helmholtz equation. Numerical experiments show that these two multilevel neural solvers effectively overcome the drawback of single-level methods, providing better robustness and computational efficiency.
\end{abstract}




\end{frontmatter}



\section{Introduction}\label{sec:intro}
Partial differential equations (PDEs) serve as fundamental mathematical models for describing scientific and engineering phenomena. Solving linear PDEs typically involves numerical methods, such as finite difference and finite element methods (FEM) \cite{quarteroni2008numerical}, which discretize these equations into large-scale linear systems. The solution for linear systems  present a bottleneck in scientific and engineering computations. Iterative methods \cite{saad2003iterative} are widely used to address such problems due to their low memory requirements and computational efficiency, and can be categorized into single-level and multilevel approaches. Single-level methods, such as Richardson and weighted Jacobi iterative method, have low per-step costs but often converge slowly. To improve convergence rates, momentum techniques, like Nesterov accelerated gradient (NAG) \cite{nesterov2003introductory, nesterov1983}, have been applied. However, the convergence rate is still constrained by the condition number of the coefficient matrix.
To overcome this limitation, preconditioning techniques are introduced, with multilevel methods \cite{trottenberg2000multigrid} being particularly important due to their enhanced computational efficiency and robustness. Despite this, the effective design of these methods often depends on the specific characteristics of the problem. Recently, with neural networks (NNs) being widely applied across various fields, improving the efficiency and robustness of iterative methods has become a significant research focus in scientific computing.

Neural networks have been proven to possess powerful expressive capabilities \cite{e2019barron} and can be used as universal approximators for PDE solutions. By minimizing loss functions derived from the strong or variational forms of PDEs, various neural PDE solvers have been introduced, including physics-informed neural networks \cite{lu2021physics, raissi2019physics} and the deep Ritz method \cite{yu2018deep}. While these neural PDE solvers have successfully tackled complex problems, they also face certain limitations. 
On the one hand, when PDE parameters (such as boundary conditions) change, these solvers often require retraining, making it costly to address parametric PDEs. On the other hand, achieving the same high accuracy like traditional methods remains challenging for these approaches. 
To address the first issue, neural operators for parametric PDEs have been developed to directly learn the solution map from the parameter space to the solution space, such as DeepONet \cite{lu2019deeponet} and the Fourier Neural Operator (FNO) \cite{li2020fourier}. However, the accuracy of neural operators still aligns with that of conventional neural PDE solvers. 
To tackle the second limitation, \textit{neural solvers} have emerged, leveraging deep learning techniques to accelerate traditional iterative methods, thereby enabling high-accuracy, efficient solutions for parametric PDEs.

Neural solvers can be broadly categorized into two approaches. The first involves using neural networks to predict key parameters in iterative methods. For example, neural networks have been employed to learn parameters in Chebyshev semi-iterative methods to enhance smoothing effects \cite{cui2024neural}, and to determine iterative directions in conjugate gradient methods \cite{kaneda2022deep}. In preconditioning techniques, neural networks can be used to correct incomplete LU decompositions \cite{trifonov2024learning}, learning smoothers \cite{katrutsa2020black, huang2022learning, chen2022meta}, transfer operators \cite{greenfeld2019learning, luz2020learning, kopanivcakova2024deeponet, liu_learning_2024}, and coarsening strategies \cite{taghibakhshi2021optimization, caldana2023deep, zou2023autoamg} in multigrid method. Moreover, neural networks can learn adaptive coarse grid basis functions \cite{klawonn2024learning} and optimize boundary conditions and interpolation operators \cite{taghibakhshi2022learning, knoke2023domain, taghibakhshi2023mg, kopanivcakova2024deeponet}  in domain decomposition methods.
The second approach integrates traditional iterative methods with neural networks to form deep learning-based hybrid iterative methods (DL-HIM) \cite{zhang2022hybrid, cui2022fourier, cui2024convergence, hu2024hybrid, xie2023mgcnn, rudikov2024neural}. The primary motivation behind DL-HIM is to leverage the spectral bias of neural networks, which tend to learn low-frequency functions\cite{rahaman2019spectral, hong2022activation, xu2019frequency}, while traditional iterative methods are effective at smoothing high-frequency error components. 
By alternating using neural networks and simple iterative methods, DL-HIM can efficiently eliminate error components across different frequencies. Pioneering work includes HINTS \cite{zhang2022hybrid} and the Fourier neural solver (FNS) \cite{cui2022fourier, cui2024convergence}.

This paper first considers the \( m \)-step Richardson iteration as a single-step iteration, proposing the Richardson(\( m \)) iterative method with periodic weights. It has been proven that the long-step version has a faster convergence rate \cite{grimmer_provably_2024}. We then construct a neural network to predict these weights, resulting in the Richardson(\( m \)) neural solver (Richardson(\( m \))-NS). Numerical experiments on anisotropic second-order elliptic equations show that the performance of Richardson(\( m \))-NS is comparable to that of the Chebyshev iterative method with optimal weights and Chebyshev semi-iterative method, but without the need for eigenvalue computation, leading to lower computational complexity.
Subsequently, by incorporating momentum acceleration and preconditioning techniques, we construct a series of single-level neural solvers with faster convergence. For example, when \( m=3 \), the NAGex-Richardson(\( m \))-NS with NAG method and SSOR preconditioning can achieve nearly a 20 times speedup compared to Richardson(\( m \))-NS.

To address the limitations of these single-level solvers, such as their strong dependency on PDE parameters and grid size, we integrate them as learnable smoothers within recently developed multilevel neural solvers. First, for anisotropic second-order elliptic equations, we apply them to FNS as smoothers. By utilizing an alternating optimization algorithm, we obtain a Richardson(\( m \))-FNS that is independent of PDE parameters and mildly dependent on grid size, significantly reducing iteration counts compared to the original FNS. For the Helmholtz equations with high wavenumbers, we incorporate the momentum-accelerated Richardson(\( m \))-NS as a new smoother in the Wave-ADR neural solver \cite{cui2024neural}, resulting in the NAG-Richardson(\( m \))-WANS, which has a parameter-independent convergence rate and demonstrates improved performance over the original Wave-ADR neural solver.

The remainder of this paper is organized as follows: Section \ref{sec:02} introduces the Richardson(m)-NS, along with momentum acceleration and preconditioning techniques. Section \ref{sec:03} integrates these single-level neural solvers with multilevel methods. Section \ref{sec:04} concludes the paper and outlines future research directions.

\section{Richardson(m) Neural Solver and Acceleration Variants}\label{sec:02}

Consider the following linear algebraic system
\begin{equation}
    \label{chapter2:eq:linearsys}
    \bm A \bm u = \bm f,
\end{equation}
where \(\bm A \in \mathbb{C}^{N \times N}\), \(\bm f \in \mathbb{C}^N\), and \(\bm u \in \mathbb{C}^N\) is the solution vector to be determined. The iteration format of the Richardson method for solving equation \eqref{chapter2:eq:linearsys} is given by
\begin{equation}
    \begin{aligned}
        \bm u_k & = (\bm I - \omega \bm A) \bm u_{k-1} + \omega \bm f \\ 
        & = \bm u_{k-1} + \omega (\bm f - \bm A \bm u_{k-1}),
    \end{aligned}
    \label{eq:richard}
\end{equation}
where \(\bm u_k\) and \(\bm u_{k-1}\) are the solutions at the \(k\)-th and \((k-1)\)-th iterations, respectively. \(\omega\) is an adjustable weight parameter, and \(\bm I\) is the identity matrix.

It has been proven that utilizing different weights at different steps can improve the convergence rate of \eqref{eq:richard} \cite{grimmer_provably_2024}. 
Specifically, we regard the continuous  $m$ (with $m \geq 1$) iterations of \eqref{eq:richard} (referred to as inner iterations) as one outer iteration. For each inner iteration, a separate weight is assigned, resulting in the following iteration format
\begin{equation}
    \label{chapter2:eq:Richardson(m)}
    \bm u_k = \bm T_m(\bm A, \bm \omega) \bm u_{k-1} + \bm g,
\end{equation}
where $\bm \omega = (\omega_1, \ldots, \omega_m)$,
\begin{equation}
    \label{chapter2:eq:Chebyshev(m)T}
    \bm T_m(\bm A, \bm \omega) = (\bm I - \omega_m \bm A)(\bm I - \omega_{m-1} \bm A) \cdots (\bm I - \omega_1 \bm A),
\end{equation}
\begin{equation}
    \label{chapter2:eq:Chebyshev(m)f}
    \bm g = \left( \omega_m + \sum_{j=1}^{m-1} \prod_{i=j+1}^{m} (\bm I - \omega_i \bm A) \omega_j \right) \bm f.
\end{equation}
This iteration method is termed the Richardson($m$) method since each outer iteration employs $m$ different weights.

The choice of the weight parameter $\bm \omega$ directly impacts the convergence speed of Richardson($m$). Inappropriate weights may result in slow convergence or even divergence of the iterative method. Consequently,  how to rationally select weights to ensure fast convergence of Richardson($m$)  is an  important issue.
In fact, when $\bm A$ is a symmetric positive definite matrix, this problem has been theoretically solved. The authors of \cite{golub1961chebyshev, axelsson1990algebraic} provide a method for constructing the weights $\bm \omega$ using Chebyshev optimal approximation theory \cite{mason2002chebyshev}. Here, we present the results  derived using this method. Assuming that the maximum and minimum eigenvalues of $\bm A$, $\lambda_{\max}$ and $\lambda_{\min}$, are known, the weights $\omega_i \in \bm \omega$ can be selected  as
\begin{equation}\label{eq:chebyshev-weights}
  \omega_i = \frac{2}{\lambda_{\max} + \lambda_{\min} + (\lambda_{\max} - \lambda_{\min}) x_i}, \quad i=1,\ldots,m,
\end{equation}
where
\begin{equation}\label{eq:chebyshev-weights-cos}
  x_i = \cos \frac{(2i-1)\pi}{2m}, \quad i=1,\ldots,m,
\end{equation}
is the root of the $m$-th order Chebyshev polynomial.

However, computing the minimum eigenvalue of $\bm A$ often requires a very high computational cost, which means that the expression in \eqref{eq:chebyshev-weights} is of little practical value. A compromise is to replace $\lambda_{\min}$ with $\alpha \lambda_{\max}$, where the constant factor $0 < \alpha < 1$ can be selected empirically; for example, $\alpha = 1/30$ in \cite{adams2003parallel}. The Richardson(m) method using this approach to obtain weights is referred to as the Chebyshev semi-iterative method.

Nevertheless, the Chebyshev semi-iterative method remains unsatisfactory, as it requires not only determining the maximum eigenvalue but also involves multiple attempts to select $\alpha$. Furthermore, if $\bm A$ is not a symmetric positive definite matrix, it is still difficult to provide a set of universal weight expressions for Richardson(m) from a theoretical standpoint. In light of these issues, this paper use neural networks to learn the parameters $\bm \omega$, aiming to obtain a faster-converging Richardson(m) neural solver. Before introducing the learning algorithm, we will first present several versions of Richardson(m) based on momentum acceleration.

\subsection{Momentum-Accelerated Richardson(m)}
In this section, we propose several momentum-based Richardson(m) methods. To this end, we first introduce an important optimization method: the Gradient Descent (GD) method, highlighting that GD and the Richardson method are formally identical, even though they are applied to solve different types of problems.

Unconstrained optimization problems can typically be represented as
\begin{equation}\label{eq:optimization-problem}
    \mathop{\min}\limits_{\bm \theta \in \Theta} J(\bm \theta),
\end{equation}
where $\bm \theta$ is the optimization variable and $J$ is the objective function. GD is a widely used optimization algorithm. Assuming that $J$ is a continuous function, the GD for solving problem \eqref{eq:optimization-problem} can be expressed as
\begin{equation}\label{eq:gd}
    \bm \theta_{k} = \bm \theta_{k-1} - \beta \nabla J(\bm \theta_{k-1}),
\end{equation}
where $\beta > 0$ is known as the step size or learning rate.

If $\bm A \in \mathbb{R}^{N \times N}$ is symmetric positive definite, we define
\begin{equation}\label{eq:Jv-linear-system}
    J(\bm u)=\frac{1}{2}\bm u^{\top}\bm A \bm u - \bm f^{\top} \bm u,
\end{equation}
then the system of equations \eqref{chapter2:eq:linearsys} is equivalent to
\begin{equation}\label{eq:Ju}
    \bm u = \arg\mathop{\min}\limits_{\bm v \in \mathbb{R}^{N}} J(\bm v).
\end{equation}
It is easy to see that $\nabla J(\bm u) = \bm A \bm u - \bm f$. Writing the corresponding GD for \eqref{eq:Ju}, we immediately obtain the Richardson method.

Introducing momentum into GD to accelerate its convergence speed is an important development direction for GD methods, leading to various variants. Momentum \cite{qian1999momentum} and NAG \cite{nesterov2003introductory, nesterov1983} are representative of these methods. Momentum adds an exponentially decaying moving average of the previous iterations' gradients to the current iteration's velocity $\bm v$. Building on momentum, NAG temporarily updates the parameters $\bm \theta$ before calculating the gradient, making the gradient update direction more predictive. These two methods are displayed in the first column of Table~\ref{tb:optimization-iterative-method}. It is noteworthy that in the last row of Table~\ref{tb:optimization-iterative-method}, we introduce a new NAG method referred to as NAG-extension. If we set $\tilde{\alpha} = \alpha$, NAG-extension degenerates into NAG; if we set $\tilde{\alpha} = 1$, we obtain the Regularised Update Descent \cite{Botev2016NesterovsAG}.

Inspired by the improvements of GD through momentum and NAG, we apply these improvement ideas to accelerate Richardson(m), proposing several momentum-based versions of Richardson(m). These methods are listed in the last column of Table~\ref{tb:optimization-iterative-method}.
\begin{table}[!htbp]
  \centering
  \caption{Optimization methods and iterative methods}\label{tb:optimization-iterative-method}
  \scalebox{0.9}[0.9]{
  \begin{tabular}{c|cc}
    \hline
    Optimization methods  & \multicolumn{2}{c}{Iterative methods}   \\
    \hline
    (GD) & (Richardson) & (Richardson(m))  \\
    $\begin{cases}
        \bm v_{k} = - \beta \nabla J(\bm \theta_{k-1})
        \\
        \bm \theta_{k} = \bm \theta_{k-1} + \bm v_{k}
     \end{cases}$
    & 
    $\begin{cases}
        \bm v_{k} = \omega(\bm f-\bm A \bm u_{k-1})
        \\
        \bm u_k = \bm u_{k-1} + \bm v_k
     \end{cases}$
    &
    $\begin{cases}
        \bm v^i_k = \omega_i(\bm f-\bm A \bm u^{i-1}_{k})
        \\
        \bm u^i_k = \bm u^{i-1}_{k} + \bm v^i_k
     \end{cases}$
    \\
     &  &  $i = 1,\ldots,m,~ \bm u^0_{k+1} = \bm u^m_k$  \\
    \hline
    (Momentum) & (Momentum) & (MOM-Richardson(m))  \\
    $\begin{cases}
        \bm v_k = {\color{red} \alpha \bm v_{k-1}} - \beta \nabla J(\bm \theta_{k-1})
        \\
        \bm \theta_k = \bm \theta_{k-1} + \bm v_k
     \end{cases}$
    & 
    $\begin{cases}
        \bm v_k = {\color{red} \alpha \bm v_{k-1}} + \omega(\bm f-\bm A \bm u_{k-1})
        \\
        \bm u_k = \bm u_{k-1} + \bm v_k
     \end{cases}$
    &
    $\begin{cases}
        \bm v^i_k = {\color{red} \alpha_i \bm v^{i-1}_{k}} + \omega_i(\bm f-\bm A \bm u^{i-1}_{k})
        \\
        \bm u^i_k = \bm u^{i-1}_{k} + \bm v^i_k
     \end{cases}$
    \\
     &  &  $i = 1,\ldots,m,~ \bm u^0_{k+1} = \bm u^m_k,~ \bm v^0_{k+1} = \bm v^m_k$  \\
    \hline
    (NAG) & (NAG) & (NAG-Richardson(m))  \\
    $\begin{cases}
        {\color{red} \tilde{\bm \theta}_k} = \bm \theta_{k-1} + \alpha \bm v_{k-1}
        \\
        \bm v_k = \alpha \bm v_{k-1} - \beta \nabla J({\color{red} \tilde{\bm \theta}_k})
        \\
        \bm \theta_k = \bm \theta_{k-1} + \bm v_k
     \end{cases}$
    & 
    $\begin{cases}
        {\color{red} \tilde{\bm u}_k} = \bm u_{k-1} + \alpha \bm v_{k-1}
        \\
        \bm v_k = \alpha \bm v_{k-1} + \omega(\bm f-\bm A {\color{red} \tilde{\bm u}_k})
        \\
        \bm u_k = \bm u_{k-1} + \bm v_k
     \end{cases}$
    &
    $\begin{cases}
        {\color{red} \tilde{\bm u}^i_k} = \bm u^{i-1}_{k} + \alpha_i \bm v^{i-1}_{k}
        \\
        \bm v^i_k = \alpha_i \bm v^{i-1}_{k} + \omega_i(\bm f-\bm A {\color{red} \tilde{\bm u}^i_k})
        \\
        \bm u^i_k = \bm u^{i-1}_{k} + \bm v^i_k
     \end{cases}$
    \\
     &  &  $i = 1,\ldots,m,~ \bm u^0_{k+1} = \bm u^m_k,~ \bm v^0_{k+1} = \bm v^m_k$  \\
    \hline
    (NAG-extension) & (NAG-extension) & (NAGex-Richardson(m))  \\
    $\begin{cases}
        \tilde{\bm \theta}_k = \bm \theta_{k-1} + {\color{red} \tilde{\alpha}} \bm v_{k-1}
        \\
        \bm v_k = \alpha \bm v_{k-1} - \beta \nabla J(\tilde{\bm \theta}_k)
        \\
        \bm \theta_k = \bm \theta_{k-1} + \bm v_k
     \end{cases}$
    & 
    $\begin{cases}
        \tilde{\bm u}_k = \bm u_{k-1} + {\color{red} \tilde{\alpha}} \bm v_{k-1}
        \\
        \bm v_k = \alpha \bm v_{k-1} + \omega(\bm f-\bm A \tilde{\bm u}_k)
        \\
        \bm u_k = \bm u_{k-1} + \bm v_k
     \end{cases}$
    &
    $\begin{cases}
        \tilde{\bm u}^i_k = \bm u^{i-1}_{k} + {\color{red} \tilde{\alpha}_i} \bm v^{i-1}_{k}
        \\
        \bm v^i_k = \alpha_i \bm v^{i-1}_{k} + \omega_i(\bm f-\bm A \tilde{\bm u}^i_k)
        \\
        \bm u^i_k = \bm u^{i-1}_{k} + \bm v^i_k
     \end{cases}$
    \\
     &  &  $i = 1,\ldots,m,~ \bm u^0_{k+1} = \bm u^m_k,~ \bm v^0_{k+1} = \bm v^m_k$  \\
    \hline
  \end{tabular}
  }
\end{table}

\begin{remark}
    We can understand the development process of the methods in the last column of Table~\ref{tb:optimization-iterative-method} from two perspectives. First, looking from top to bottom, we successively add standard momentum and Nesterov's momentum to Richardson(m). Secondly, looking from left to right, it can be considered as transforming a single-step iteration into a long-step iteration, where consecutive m iterations are bundled together as a single outer iteration.
\end{remark}

Preconditioning techniques are commonly used in the solution of discrete systems. Another method we have in mind to accelerate these iterative methods is to add an appropriate preconditioning operator $\bm B^{-1}$, which will provide a more general iterative framework. For instance, the MOM-Richardson(m) with preconditioning can be expressed as
\begin{equation}\label{chap2:eq:momentum-m-Binv}
    \begin{cases}
        \bm v^i_k = \alpha_i \bm v^{i-1}_{k} + \omega_i\bm B^{-1}(\bm f-\bm A \bm u^{i-1}_{k})
        \\ 
        \bm u^i_k = \bm u^{i-1}_{k} + \bm v^i_k
    \end{cases}
    \quad i = 1,\ldots,m,~ \bm u^0_{k+1} = \bm u^m_k,~ \bm v^0_{k+1} = \bm v^m_k.
\end{equation}
Similarly, the versions of other iterative methods in Table~\ref{tb:optimization-iterative-method} with preconditioning can be derived. Setting $\bm B^{-1}=\bm I$ will yield the iterative methods presented in the table.

Next, we introduce some specific preconditioners. It is well known that common stationary iterative methods can be expressed as Richardson methods with preconditioning
\begin{equation}\label{chap2:eq:pre_Richardson}
    \bm u_{k} = \bm u_{k-1}+ \omega \bm B^{-1}(\bm f-\bm A\bm u_{k-1}).
\end{equation}
For example, when $\omega=1$, the corresponding preconditioners $\bm B^{-1}$ for weighted Jacobi, Gauss-Seidel, SOR, and SSOR iterative methods are $\alpha\bm D^{-1}$, $(\bm D-\bm L)^{-1}$, $\alpha(\bm D-\alpha \bm L)^{-1}$, and $\alpha(2-\alpha)(\bm D-\alpha \bm U)^{-1}\bm D(\bm D-\alpha \bm L)^{-1}$, respectively, where $\bm D,-\bm U,-\bm L$ are the diagonal, upper, and lower triangular matrices of $\bm A$, and $\alpha$ is the relaxation factor of the iterative method. 
We can apply these preconditioners to the iterative methods in Table~\ref{tb:optimization-iterative-method}. Additionally, we can also construct some problem-specific preconditioners by integrating the physical background of the linear system to be solved.

\begin{remark}
    After adding a certain stationary iterative method's preconditioner to the iterative methods in Table~\ref{tb:optimization-iterative-method}, we will obtain a specific new iterative method. We can understand the generation process of the new iterative method from many angles. For example, taking the preconditioner in equation \eqref{chap2:eq:momentum-m-Binv} to be the SOR corresponding preconditioner $\bm B^{-1}=\alpha(\bm D-\alpha \bm L)^{-1}$, and naming the new iterative method as MOM-Richardson(m)-SOR, one perspective for understanding this iterative method could be that we first introduce momentum in SOR and then transform the iterative method into a long-step version.
\end{remark}

So far, we have introduced various momentum methods and preconditioners to accelerate Richardson($m$) and obtained a series of iterative methods. These newly proposed methods share a common characteristic: they contain a large number of weights that need to be determined. Specifying the parameters of these methods is an urgent problem that requires addressing. Next, we utilize neural networks to learn these weight parameters.

\subsection{A Neural Solver Framework for Iterative Method}

Now, we present the algorithm framework for using neural networks to learn the weights of the iterative methods in Table~\ref{tb:optimization-iterative-method}. In the following sections, we will refer to the iterative method that uses NN to obtain weight parameters as the corresponding Neural Solver (NS), such as Richardson(m)-NS. For ease of description, we will denote the weights involved in these iterative methods uniformly as $\bm \omega = (\omega_1,\ldots,\omega_m, \alpha_1,\ldots,\alpha_m, \tilde{\alpha}_1,\ldots,\tilde{\alpha}_m)$, where $m\geq 1$. 

Neural networks are widely used for learning various types of linear and nonlinear mappings due to their powerful approximation capabilities.
Before applying neural networks, we briefly discuss the input  parameters for learning the mapping to the weights $\bm \omega$. Generally, $\bm \omega$ is related to the matrix A, suggesting that a functional relationship can be learned from A to $\bm \omega$.Currently, we are focused on solving parameterized PDEs, with the parameters in these PDEs uniformly represented by $\bm \mu$. Once the discrete format of the PDE is selected, $\bm \mu$ corresponds uniquely to the coefficient matrix $\bm A$ of the discrete system. Therefore, we can also use NNs to learn the functional relationship from $\bm \mu$ to $\bm \omega$.

In this paper, the NN used to learn the weights $\bm \omega$ is named Meta$(\bm \theta;\bm \mu)$, where $\bm \theta$ represents the parameters of the NN. No specific type of network model is designated for Meta; in fact, we believe that classic NNs in Deep Learning, such as fully-connected NNs (FNNs) and convolutional NNs (CNNs), can all be attempted. Users can design a specific network structure based on the characteristics of the input and output data and other requirements.

Typically, there are two loss functions available, namely the relative error and the relative residual, expressed as
\begin{align}
    \mathcal{L}_1 &= \frac{1}{N_b}\sum^{N_b}_{i=1}\frac{\|\bm u_i - \bm u^K_i\|}{\|\bm u_i\|},  \label{eq:loss-fun-relative-err} \\ 
    \mathcal{L}_2 &= \frac{1}{N_b}\sum^{N_b}_{i=1}\frac{\|\bm f_i - \bm A_i \bm u^K_i\|}{\|\bm f_i\|},  \label{eq:loss-fun-relative-res}
\end{align}
where $N_b$ is the batch size, $\bm A_i, \bm f_i$ is the $i$-th discrete system in a batch of data, $\bm u_i$ is the corresponding true solution, and $\bm u^K_i$ denotes the numerical solution obtained after $K$ iterations starting from a given initial guess $\bm u^0_i$. Generally speaking, the true solution $\bm u_i$ is unknown. If we want to use $\mathcal{L}_1$ as the loss function, we can obtain $\bm u_i$ by solving the linear system using direct methods. However, this approach incurs a huge computational cost. Therefore, we recommend using $\mathcal{L}_2$.

Here, we present the training algorithm for the NS corresponding to the iterative methods in Table~\ref{tb:optimization-iterative-method}.

\begin{algorithm}[!htbp]
    \caption{Training algorithm for NS}
    \label{alg:NS-framework}
    \KwData{PDE parameters $\{\bm \mu_i\}^{N_{\text{train}}}_{i=1}$ and corresponding discrete systems $\{\bm A_i, \bm f_i\}^{N_{\text{train}}}_{i=1}$, solution $\{\bm u_i\}^{N_{\text{train}}}_{i=1}$(if using $\mathcal{L}_1$)}
    \KwIn{Meta$(\bm \theta;\bm \mu)$, $N_{\text{epochs}}$, $N_{\text{b}}$, $K$, m}
    \For{$l=1,\ldots,N_{\text{epochs}}$} {
        \For{$j=1,\ldots,\lfloor N_{\text{train}}/N_{\text{b}}\rfloor$} {
            Sample $\{\bm \mu_i\}^{N_{\text{b}}}_{i=1} \subset \{\bm \mu_i\}^{N_{\text{train}}}_{i=1}$, $\{\bm A_i, \bm f_i\}^{N_{\text{b}}}_{i=1} \subset \{\bm A_i, \bm f_i\}^{N_{\text{train}}}_{i=1}$, $\{\bm u_i\}^{N_{\text{b}}}_{i=1} \subset \{\bm u_i\}^{N_{\text{train}}}_{i=1}$ \;
            Obtain $\bm \omega=\text{Meta}(\bm \theta;\bm \mu)$ using current $\bm \theta$ \;
            \For{$i=1,\ldots,N_{\text{b}}$} {
                Obtain the initial guess $\bm u^0_i$ \;
                \For{$k=1,\ldots,K$}{
                    Compute $\bm u^K_i$ based on the iterative method in Table~\ref{tb:optimization-iterative-method} \;
                }
                Calculate the loss function $\mathcal{L}$ as follows:
                \begin{equation}
                    \mathcal{L} = 
                    \begin{cases}
                        \mathcal{L}_1 & \text{if using } \mathcal{L}_1 \\
                        \mathcal{L}_2 & \text{if using } \mathcal{L}_2
                    \end{cases}
                \end{equation}
                Backpropagate the error to update $\bm \theta$;
            }
        }
    }
    \KwOut{Learned NS}
\end{algorithm}

After training, the steps to evaluate or test the neural solver are straightforward. First, input the test data, $\{\bm \mu_i\}^{N_{\text{test}}}_{i=1}$, into the learned Meta$(\bm \theta; \bm \mu)$ to obtain the weights $\{\bm \omega_i\}^{N{\text{test}}}_{i=1}$, and then the corresponding iterative method is executed until a certain acceptable tolerance or the maximum number of iterations is reached. We will not list the detailed algorithm again.

\subsection{Numerical Experiments}\label{sec:sec-2-ne}

In this section, we evaluate the performance of the proposed neural solvers through experiments, using an anisotropic second-order elliptic problem as an example.

Consider the following two-dimensional diffusion equation with anisotropic coefficients:
\begin{equation}\label{eq:anisotropy-poisson}
    \begin{cases}
    \begin{aligned}
      -\nabla\cdot(C\nabla u) &= f, \quad \text{in}~ \Omega, \\ 
      u &= 0, \quad \text{on}~ \partial\Omega,
    \end{aligned} 
    \end{cases}
\end{equation}
where
\begin{equation*}
  C=C(\varepsilon,\theta)=
    \begin{pmatrix}
    \cos \theta & -\sin \theta \\ 
    \sin \theta & \cos \theta
    \end{pmatrix}
    \begin{pmatrix}
    1 & 0 \\ 
    0 & \varepsilon
    \end{pmatrix}
    \begin{pmatrix}
    \cos \theta & \sin \theta \\ 
    -\sin \theta & \cos \theta
    \end{pmatrix},
\end{equation*}
$0<\varepsilon<1$ is the anisotropic strength, $\theta \in[0,\pi]$ is the anisotropic direction, and $\Omega=(0,1)^2$. 

We create a uniform grid partition of $\Omega$ into $N\times N$ and use bilinear elements to discretize ~\eqref{eq:anisotropy-poisson}~. After handling the boundaries, we obtain the discrete system shown in \eqref{chapter2:eq:linearsys}.

Before training the model, we need to generate some training and testing data.
In this example, the PDE parameters are $\bm \mu = (\varepsilon,~ \theta)$. In the training dataset, the sample size is $N_{\text{train}}=10000$, and we use a uniform random number generator to obtain $\bm \mu$, where $\varepsilon \sim \mathcal{U}(10^{-6},1),~ \theta \sim \mathcal{U}(0,\pi)$. The grid partition size is fixed at $N=64$. With $\bm \mu$ obtained, we can easily derive the corresponding matrix $\bm A$ using the bilinear element format. Next, we generate the right-hand side vector using a normal distribution, $\bm f \sim \mathcal{N}(\bm 0, \bm I)$. In the testing phase, we select some specific $\bm \mu$ and $N$ according to the experimental purpose, then generate the corresponding $\bm A$ and $\bm f$ in the same way.
It is worth noting that, to enhance the network's sensitivity to the value of $\varepsilon$, the input to the network is actually $\lg\varepsilon$ instead of $\varepsilon$.

\begin{figure}[!htbp]
    \centering
    \includegraphics[width=0.6\textwidth]{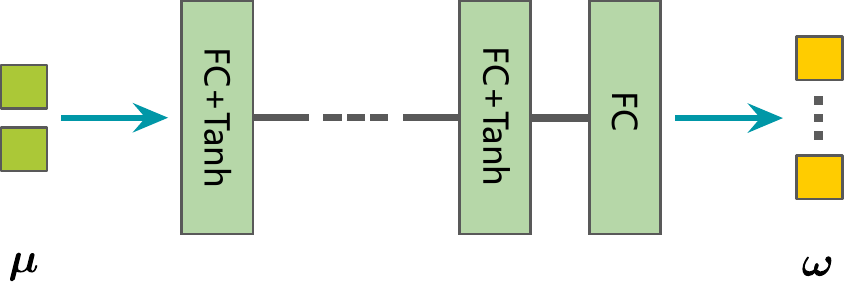}
    \caption{A schematic of a Meta$(\bm \theta;\bm \mu)$ network, i.e., an FNN.}
    \label{fig:meta-fnn}
\end{figure}

We select the loss function $\mathcal{L}_2$ as shown in \eqref{eq:loss-fun-relative-res}. Since both the input and output data of the network, $\bm \mu$ and $\bm \omega$, are one-dimensional, we design Meta$(\bm \theta;\bm \mu)$ as a simple FNN, as shown in Figure~\ref{fig:meta-fnn}.
During the validation or testing process, we use the relative residual $\|\bm f - \bm A \bm u_k\|/\|\bm f\|$ to measure the approximation accuracy of the numerical solution to the true solution. For solvers with converging iterative processes, the iteration will stop if the relative residual is less than $10^{-6}$.

Based on the aforementioned training data, we trained the NS corresponding to several iterative methods listed in Table~\ref{tb:optimization-iterative-method} using Algorithm~\ref{alg:NS-framework}. After training, we evaluated the performance of these solvers on some discrete systems in the test set and recorded the number of iterations required to reach the threshold in Tables~\ref{tb:Richardson-m},~\ref{tb:richa-different-parameters}, and~\ref{tb:richa-different-scales}. All these solvers started the iterations with a zero initial guess.
All numerical experiments in this paper were conducted on an Nvidia A100-SXM4-80GB GPU.

From Table~\ref{tb:Richardson-m}, it can be seen that the Richardson(m)-NS achieved a number of iterations close to the Chebyshev iteration with the same $m$, and it outperformed the Chebyshev semi-iteration method. Although this phenomenon was observed for $m=3$, it holds for other values of $m$ as well. Observing the results for Richardson(m)-NS at $m=1,3,7,15$, we can see that increasing $m$ significantly improves the convergence speed. In fact, the same pattern is observed for other solvers listed in Table~\ref{tb:optimization-iterative-method}. Due to space limitations, we do not list the results. Comparing the results of Richardson(m)-NS, MOM-Richardson(3)-NS and NAG-Richardson(3)-NS, we can see that the introduction of both standard momentum and Nesterov's momentum in Richardson(m) greatly improves the convergence rate, and the latter has a better acceleration effect. NAGex-Richardson(3)-NS achieved better results compared to NAG-Richardson(3)-NS, which indicates that our correction to NAG is beneficial. On this basis, the convergence speed of the solver is further significantly enhanced by incorporating the SSOR preconditioner $\bm B^{-1}_{SSOR}$.

Number of iterations of solvers in table ~\ref{tb:Richardson-m}~ increases as $\varepsilon$ decreases. In table ~\ref{tb:richa-different-parameters}~,  iteration counts fluctuates unsteadily with the increase of $\theta$. Iterations in table \ref{tb:richa-different-scales} increases with the increase of discrete scale $N$.
These phenomena show that the convergence rate of these single-level solvers is greatly affected by PDE parameters and discrete scale.

\begin{table}[!htbp]
  \centering
  \caption{Number of iterations for various iterative methods solving the discrete system of anisotropic problems ($\theta=0$, $N=64$).}\label{tb:Richardson-m}
  \arrayrulewidth=0.8pt
  \begin{tabular}{l|ccccc}
    \hline
    \multirow{2}*{Solvers} & \multicolumn{5}{c}{$\epsilon$} \\ 
    \cline{2-6}
     &  1 & $10^{-2}$ & $10^{-4}$ & $10^{-6}$ & $10^{-8}$ \\ 
    \hline
    Richardson(1)-NS & 3999 & 7534 & 18179 & 27257 & 27721  \\ 
    Chebyshev iteration (m=3 in \eqref{eq:chebyshev-weights})  & 587 & 1105 & 2204 & 3273 & 3283  \\ 
    Chebyshev semi-iteration (m=3, $\alpha=1/30$)  & 585 & 1120 & 2708 & 4061 & 4129  \\ 
    Richardson(3)-NS  & 568 & 947 & 2223 & 3315 & 3366  \\ 
    Richardson(7)-NS  & 105 & 191 & 457 & 685 & 696  \\ 
    Richardson(15)-NS  & 29 & 50 & 122 & 184 & 187  \\ 
    \hline
    MOM-Richardson(3)-NS  & 137 & 181 & 405 & 629 & 648  \\ 
    NAG-Richardson(3)-NS  & 102 & 141 & 344 & 531 & 540  \\ 
    NAGex-Richardson(3)-NS  & 95 & 140 & 285 & 467 & 512  \\ 
    NAGex-Richardson(3)-NS (with $\bm B^{-1}_{SSOR}$) & 20 & 19 & 29 & 43 & 38  \\ 
    \hline
  \end{tabular}
\end{table}

%

\begin{table}[!htbp]
  \centering
  \caption{Number of iterations for NAGex-Richardson(3)-NS (with $\bm B^{-1}_{SSOR}$) solving discrete systems corresponding to different anisotropic directions $\theta$ ($\epsilon=10^{-6}$, $N=64$).}\label{tb:richa-different-parameters}
  \arrayrulewidth=0.8pt
  \begin{tabular}{cccccccc}
    \toprule[1pt]
    $\theta$ & 0 & $\pi/6$ & $\pi/3$ & $\pi/2$ & $2\pi/3$ & $5\pi/6$ & $\pi$  \\
    \midrule
    Iterations & 43 & 27 & 28 & 50 & 28 & 28 & 30  \\
    \bottomrule[1pt]
  \end{tabular}
\end{table}

\begin{table}[!htbp]
  \centering
  \caption{Number of iterations for NAGex-Richardson(3)-NS (with $\bm B^{-1}_{SSOR}$) solving discrete systems with different grid sizes $N$ ($\epsilon=10^{-6}$, $\theta=\pi/10$).}\label{tb:richa-different-scales}
  \arrayrulewidth=0.8pt
  \begin{tabular}{cccccc}
    \toprule[1pt]
    $N$ & 32 & 64 & 128 & 256 & 512 \\
    \midrule
    Iterations  & 26 & 27 & 57 & 250 & 977  \\
    \bottomrule[1pt]
  \end{tabular}
\end{table}

\section{Application in Multilevel Solvers}\label{sec:03}
In this section, we apply our proposed single-level neural solvers to recently developed multilevel methods to enhance their performance. First, we integrate them with the FNS for solving anisotropic diffusion equations. Next, we combine them with the Wave-ADR neural solver (WANS) to address Helmholtz equations with high-wavenumber.

\subsection{Combination with FNS}\label{sec:03-1}
The iterative format of FNS is 
\begin{align}
    \cu{u}^{(k+\frac{1}{2})} &= \cu{u}^{(k)} + \cu{B} (\cu{f} - \cu{A u}^{(k)})\quad \text{(smoothing iteration, repeat } M \text{ times)},  \label{eq:dl-him-b} \\
    \bm u^{(k+1)} &= \bm u^{(k+\frac{1}{2})} + \mathcal{H}(\bm f - \bm A \bm u^{(k+\frac{1}{2})}) \quad  \text{ (neural iteration)}.\label{eq:dl-him-h}
\end{align}
Here, the smoother $\cu{B}$ represents the weighted Jacobi preconditioner, and the smoothing iteration \eqref{eq:dl-him-b} is used to reduce certain error components, primarily high-frequency errors. The neural iteration \eqref{eq:dl-him-h} introduces $\mathcal{H}$, a neural operator designed to address error components that $\bm{B}$ cannot effectively eliminate, its formulation is given by
\begin{equation*}
    \mathcal{H} = \mathcal{F} \mathcal{T} \tilde{\Lambda} \mathcal{T}^* \mathcal{F}^{-1},
\end{equation*}
where $\mathcal{F}$ is the Fourier matrix, $\tilde{\Lambda}$ is a diagonal matrix produced by the Meta-$\tilde{\Lambda}$ network to approximate the inverses of the eigenvalues of $\bm{A}$, and $\mathcal{T}$ is a sparse transition matrix with kernel generated by the Meta-$T$ network, such that $\mathcal{F} \mathcal{T}$ can approximate the eigenvector matrix of $\bm{A}$.

\begin{figure}[!htbp]
  \centering
  \includegraphics[width=0.9\textwidth]{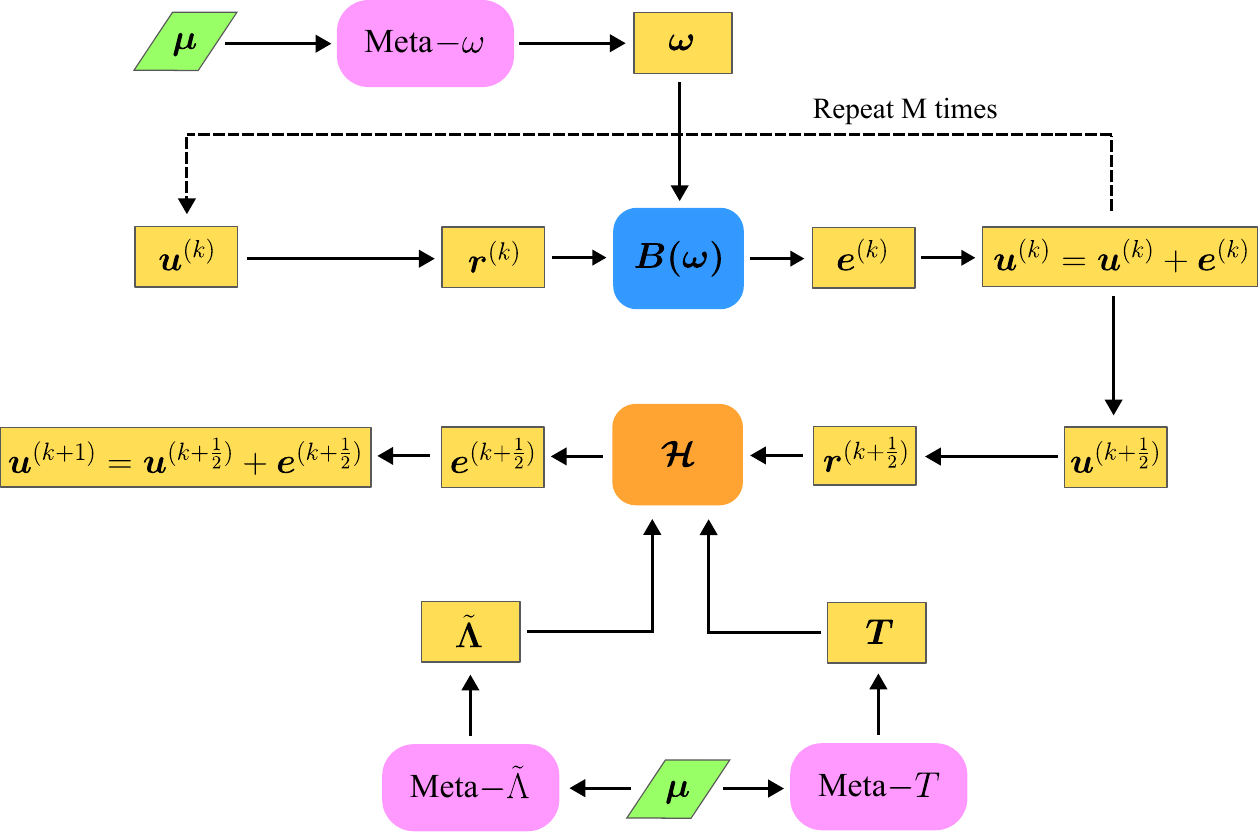}
  \caption{Computational flow of Richardson(m)-FNS.}
  \label{fig:rm-fns}
\end{figure}

\begin{figure}[!htbp]
  \centering
  \includegraphics[width=0.9\textwidth]{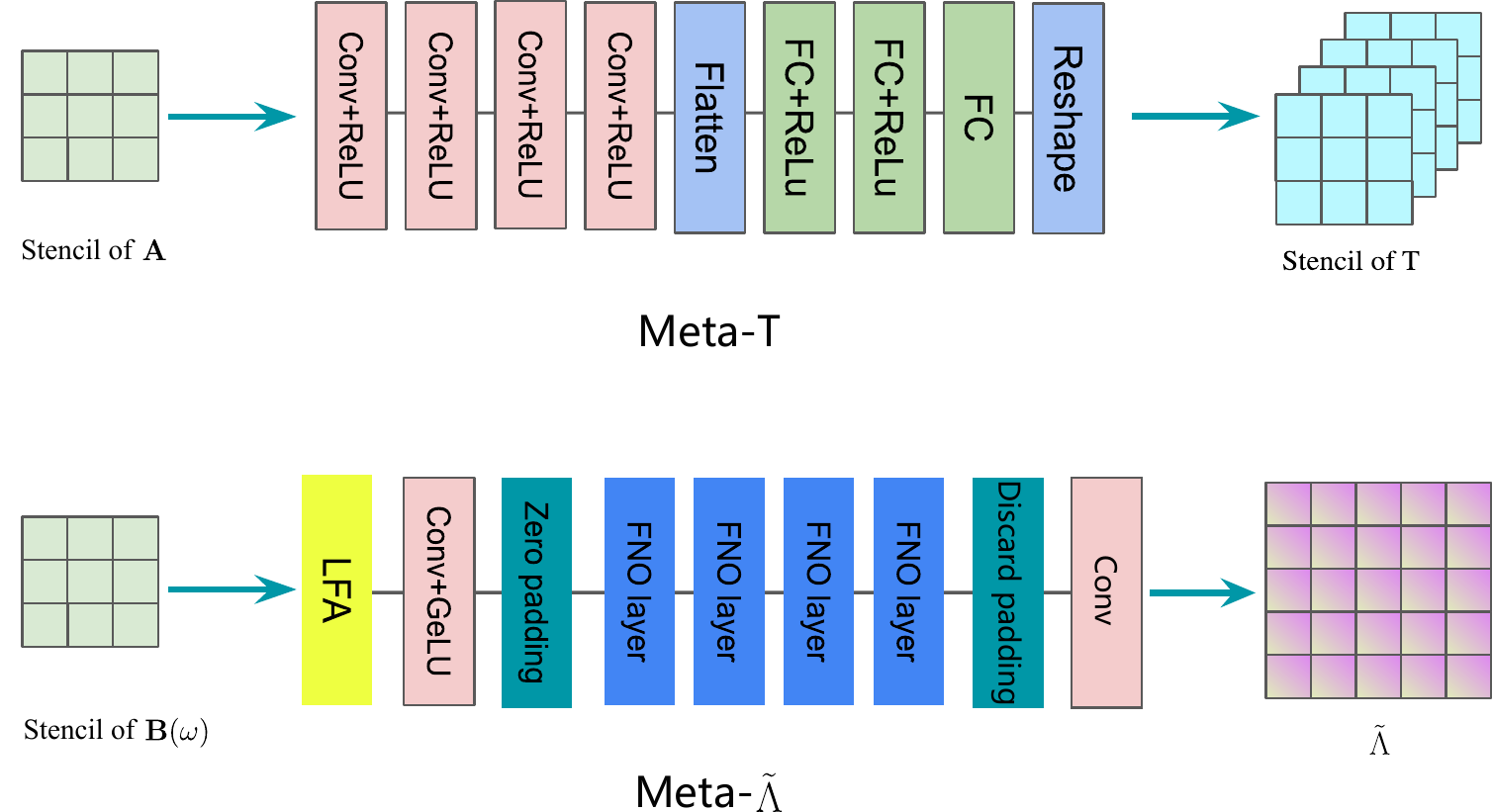}
  \caption{Network architecture of Meta-$T$ and Meta-$\tilde{\Lambda}$.}
  \label{fig:meta-2net}
\end{figure}

Since $\bm{B}$ is a fixed smoother, its smoothing effect on anisotropic problems is limited, which complicates the learning for $\mathcal{H}$. 
To address this limitation, we replace $\bm{B}$ with Richardson(m)-NS, which enhances the smoothing effect. The resulting solver is termed Richardson(m)-FNS, and its computational flowchart is shown in Figure~\ref{fig:rm-fns}. 
There are three meta networks used to learn the parameters required for Richardson(m)-FNS. Meta-\(\omega\) is a FNN depicted in Figure~\ref{fig:meta-fnn}, while Meta-\(\tilde{\Lambda}\) and Meta-T are implemented as a FNO \cite{li2020fourier} and a CNN \cite{cui2024convergence}, respectively, with the network architecture illustrated in Figure~\ref{fig:meta-2net}.

\subsubsection{Data and Training Algorithm}  
We  use the anisotropic diffusion equation \eqref{eq:anisotropy-poisson} to assess the performance of Richardson(m)-FNS. The train and test data in the following experiments remain identical to those previously used.

We first give an algorithm for training the network parameters. Richardson(m)-FNS includes three meta networks: Meta-$\omega$, Meta-$\tilde{\Lambda}$, and Meta-$T$, with a total parameter counts approaching 20 million. Based on our experience, simultaneously training these networks is challenging for achieving optimal performance, so we propose a new training Algorithm \ref{alg:alternating-optimization} utilizing an alternating optimization strategy.
\begin{algorithm}[!htbp]
  \caption{Alternating Optimization Training Algorithm}
  \label{alg:alternating-optimization}
  \KwIn{Number of alternating training cycles $M$, Number of epochs between cycles $N_{\text{interval}}$}
  Freeze Meta-$\omega$, and use a fixed $\bm{\omega}^{(0)}=(\omega^{(0)}_1,\ldots,\omega^{(0)}_m)$ to replace the output of Meta-$\omega$\;
  Perform $N_{\text{interval}}$ training epochs to obtain the network parameters $\Theta^{(0)}_{\text{Meta-}\tilde{\Lambda}}$ and $\Theta^{(0)}_{\text{Meta-}T}$, then activate Meta-$\omega$\;
  \For{$i=1,\ldots,M$}{
      Transfer and freeze $\Theta^{(i-1)}_{\text{Meta-}\tilde{\Lambda}}$ and $\Theta^{(i-1)}_{\text{Meta-}T}$, and perform $N_{\text{interval}}$ training epochs to obtain $\Theta^{(i)}_{\text{Meta-}\omega}$\;
      Transfer and freeze $\Theta^{(i)}_{\text{Meta-}\omega}$, and perform $N_{\text{interval}}$ training epochs to obtain $\Theta^{(i)}_{\text{Meta-}\tilde{\Lambda}}$ and $\Theta^{(i)}_{\text{Meta-}T}$\;
      \If{the stopping condition of training is met}{
          BREAK\;
      }
  }
  \KwOut{Learned Richardson(m)-FNS}
\end{algorithm}

\begin{remark}
Some explanations on Algorithm~\ref{alg:alternating-optimization}.
\begin{enumerate}
    \item The selection of $\bm{\omega}^{(0)}$ depends on the PDE. We set $\bm{\omega}^{(0)}=\bm{0.5}$ in the following experiments.
    \item The selection of $N_{\text{interval}}$ need not obey strict criteria; it can be determined based on the decrease in the training loss. For example, choose the number of training epochs until the loss function no longer shows significant decay.
    \item During the training process, we found that after the loss decreased to a certain level, the performance on the test set actually worsened as the loss continued to decrease. Therefore, an early stopping strategy is useful. 
\end{enumerate}
\end{remark}

Next, we test the performance of the trained Richardson(m)-FNS when solving the discretized system from the anisotropic equation. The right-hand sides are randomly generated following the standard normal distribution. In the following Tables, $N$ denotes the grid size along each direction. The stopping criterion is that the norm of the relative residual is less than $10^{-6}$.

Table~\ref{tb:rm-fns-different-m} shows the number of iterations required for Richardson(m)-FNS with different $m$ to achieve convergence when solving the discretized system from the anisotropic equation with different $\epsilon$. It is important to note that all smoothers use the same smoothing iteration count $M=10$. For instance, when $m=2$, the Richardson(m)-NS is executed for 5 iterations. As can be observed, increasing $m$ significantly improves the convergence speed. Compared to the single-level solver Richardson(m)-NS in Table~\ref{tb:Richardson-m}, the Richardson(m)-FNS exhibits a much faster convergence rate.
\begin{table}[!htbp]
  \centering
  \caption{Number of iterations for Richardson(m)-FNS to solve the linear systems corresponding to different m ($\theta=0$, $N=64$).}\label{tb:rm-fns-different-m}
  \arrayrulewidth=0.8pt
  \begin{tabular}{l|ccccc}
    \hline
    \multirow{2}{*}{Models} & \multicolumn{5}{c}{$\epsilon$}  \\
    \cline{2-6}
     &  1 & $10^{-2}$ & $10^{-4}$ & $10^{-6}$ & $10^{-8}$  \\
    \hline
    Richardson(1)-FNS  &  8 & 14 & 15 & 23 & 22  \\
    Richardson(2)-FNS  &  7 & 11 & 13 & 22 & 17  \\
    Richardson(10)-FNS & 4 & 6 & 8 & 9 & 8  \\
    \hline
  \end{tabular}
\end{table}

Tables~\ref{tb:fns-r10-fixedt-e}, \ref{tb:fns-r10-t-fixede} and \ref{tb:fns-r10-N} shows the iteration counts of Richardson(10)-FNS when solving the anisotropic equation with different parameters \(\epsilon\), \(\theta\), and grid sizes \(N\). It is evident that Richardson(10)-FNS significantly enhances the convergence rate compared to the original FNS.

\begin{table}[!htbp]
  \centering
  \caption{Iteration counts of the two FNS to solve the linear systems corresponding to different strengths $\epsilon$, where "$-$" indicates divergence ($\theta=\pi/10$, $N=64$).}\label{tb:fns-r10-fixedt-e}
  \arrayrulewidth=0.8pt
  \begin{tabular}{cccccc}
    \toprule[1pt]
    \(\epsilon\) & 1 & \(10^{-2}\) & \(10^{-4}\) & \(10^{-6}\) & \(10^{-8}\)  \\
    \midrule
    FNS  & $-$ & 13 & 16 & 14 & 13  \\
    Richardson(10)-FNS & 4 & 6 & 6 & 6 & 6  \\
    \bottomrule[1pt]
  \end{tabular}
\end{table}

\begin{table}[!htbp]
  \centering
  \caption{Iteration counts of the two FNS to solve the linear systems corresponding to different directions $\theta$ ($\epsilon=10^{-6}$, $N=64$).}\label{tb:fns-r10-t-fixede}
  \arrayrulewidth=0.8pt
  \begin{tabular}{cccccccc}
    \toprule[1pt]
    \(\theta\) & 0 & \(\pi/6\) & \(\pi/3\) & \(\pi/2\) & \(2\pi/3\) & \(5\pi/6\) & \(\pi\)  \\
    \midrule
    FNS   & 51 & 16 & 15 & 57 & 13 & 15 & 47  \\
    Richardson(10)-FNS & 9 & 6 & 6 & 10 & 6 & 5 & 9  \\
    \bottomrule[1pt]
  \end{tabular}
\end{table}
\begin{table}[!htbp]
  \centering
  \caption{Iteration counts of the two FNS to solve the linear systems corresponding to different discrete sizes $N$ ($\epsilon=10^{-6}$, $\theta=\pi/10$).}\label{tb:fns-r10-N}
  \arrayrulewidth=0.8pt
  \begin{tabular}{l|ccccc}
    \hline
    \multirow{2}*{Models} & \multicolumn{5}{c}{$N$}  \\
    \cline{2-6}
     & 32 & 64 & 128 & 256 & 512  \\
    \hline
    FNS & 14 & 16 & 16 & 43 & 178  \\
    Richardson(10)-FNS & 5 & 6 & 9 & 12 & 83  \\
    \hline
  \end{tabular}
\end{table}

\subsection{Combination with Wave-ADR}
Consider the following 2D Helmholtz equation
\begin{equation}
  \begin{cases}
  \begin{aligned}
    -\Delta u-k^2u+i\gamma\frac{\omega}{c^{2}}u &=f,\quad \text{in} \,\, \Omega,  \\
    u &= 0,  \quad \text{on}\,\,  \partial\Omega,
  \end{aligned}
  \end{cases}
\end{equation}
where $f(\boldsymbol{x})$ is the source term, $k=\omega/c(\boldsymbol{x})$ is the wavenumber, $c(\boldsymbol{x})$ is the wave speed, $\omega$ is the angular frequency, and $\gamma$ is the damping mask, which equals zero inside $\Omega$ and increases from zero to $\omega$ within the sponge layer near the boundary. The thickness of this layer is typically one wavelength\,\cite{treister2019multigrid}.

Using the second-order central finite difference method on a uniform grid, the discrete stencil at the interior nodes is
\begin{equation}\label{eq:Helmholtz-Sponge-stencil}
  \frac{1}{h^2}\left[\begin{array}{ccc}
  0 & -1 & 0 \\
  -1 & 4-k_{i,j}^2 h^2 & -1 \\
  0 & -1 & 0
  \end{array}\right].
\end{equation}
where $h=1/N$ in both the $x-$ and $y-$ directionss. 
According to the Shannon sampling principle, at least $10$ grid points per wavelength are needed, resulting in a large-scale linear system. Given the indefinite nature of the coefficient matrix, a single-level neural solver alone may cause divergence. To ensure convergence, we incorporate it with a multigrid-based neural solver, the Wave-ADR neural solver (WANS) \cite{cui2024neural}. Since standard multigrid methods often fail to eliminate errors with angular frequencies near the wavenumber, WANS introduces the ADR correction to handle these challenging error components.

\begin{figure}[!htbp]
  \centering
  \includegraphics[width=0.9\textwidth]{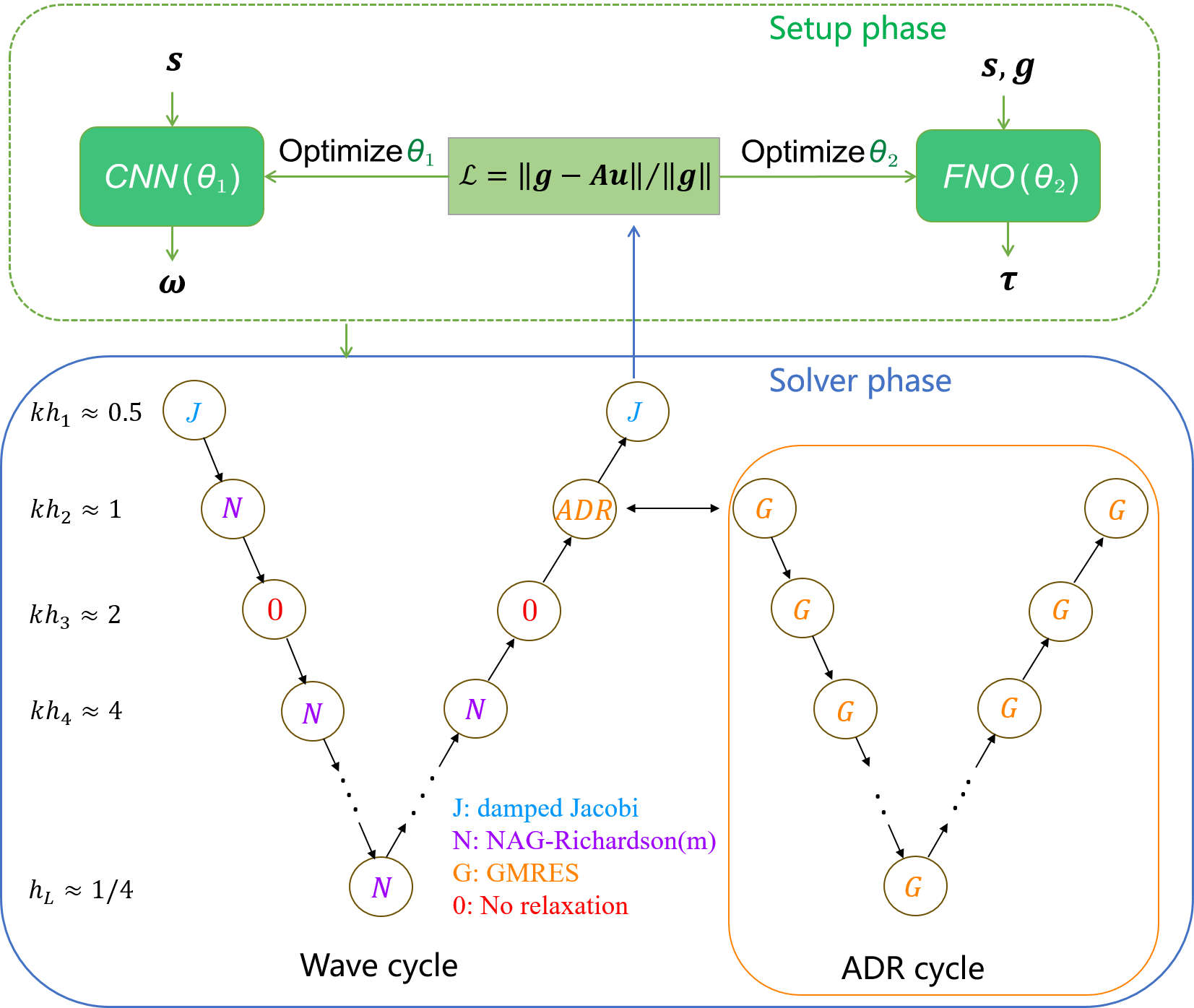}
  \caption{Computational flow of NAG-Richardson(m)-WANS.}
  \label{fig:wave-adr}
\end{figure}

In the original WANS, the Chebyshev semi-iterative method applied on the normal equations as smoother, requiring learning parameter $\alpha$ and estimating the maximum eigenvalue via the power method,  which is computationally costly for large-scale systems. To mitigate this, we use our proposed single-level neural solvers as smoothers. Figure \ref{fig:wave-adr} illustrates the computational flow of WANS with NAG-Richardson(m) as the smoother (termed NAG-Richardson(m)-WANS).

\subsubsection{Data and Experiments}
We use the same training and testing data as in the original WANS \cite{cui2024neural}. We replace the Chebyshev semi-iterative method in WANS with Richardson(m), MOM-Richardson(m), and NAG-Richardson(m). It is worth to note that we no longer rely on the network to learn the relaxation parameter $\cu{\omega}$ of the smoother directly; instead, we focus on learning $\log\cu{\omega}$. This change is based on our observation of significant differences in the parameters required for solving  linear systems at various wavenumbers (corresponding to different grid sizes), making it easier for the network to learn this new parameter. 

Finally, all methods listed in Table~\ref{tb:gmres-Wave-ADR-smoothers} are used as preconditioners for FGMRES, with a restart step of 20. We evaluate the number of iterations required by FGMRES to solve linear systems of varying sizes. It is evident that our proposed single-level neural solvers, such as NAG-Richardson(m), outperform the Chebyshev semi-iteration when used as the smoother in WANS.
\begin{table}[!htbp]
  \centering
  \caption{Iteration counts of FGMRES preconditioned by WANS with different smoothers. For all smoothers, the number of smoothing iterations is set to $m=5$.}\label{tb:gmres-Wave-ADR-smoothers}
  \arrayrulewidth=0.8pt
  \begin{tabular}{lcccccc}
    \toprule[1pt]
    $N$ & 128 & 256 & 512 & 1024 & 2048 & 4096  \\
    $\omega/2\pi$ & 10 & 20 & 40 & 80 & 160 & 320  \\
    \midrule
    WANS\cite{cui2024neural}  & 8.8 & 15.5  & 28.8  & 65.1  & 137.2  & 344.3  \\
    Richardson(5)-WANS  & 9 & 14.7 & 27.7 & 58.7 & 134.5 & 350  \\
    MOM-Richardson(5)-WANS  & 9 & 14.3  & 29  & 61.4  & 140.2 & 350.2  \\
    NAG-Richardson(5)-WANS  & 8.8 & 14.7 & 27.9 & 63.3 & 134.9  & 336.6  \\
    \bottomrule[1pt]
  \end{tabular}
\end{table}

\section{Conclusions and Future Work}\label{sec:04}


In this paper, we first design a Richardson(\( m \)) neural solver with periodic, learnable weights. We then apply momentum acceleration and preprocessing techniques to design several single-level neural solvers with improved convergence rates. Furthermore, by integrating these solvers with the recently developed FNS and WANS, we develop multilevel solvers, Richardson(\( m \))-FNS and NAG-Richardson(\( m \))-WANS, which demonstrate superior  performance compared to existing multilevel solvers when solving  discrete systems of anisotropic second-order elliptic and Helmholtz equations.
In the near future, we plan to extend these neural solvers to tackle more complex problems, such as the advection-diffusion-reaction equation.

\section*{Acknowledgments}
This work is funded by the National Natural Science Foundation of China (12371373), Science Challenge Project (Grant No. TZ2024009), Key Project of Hunan Provincial Department of Education(Grant No. 22A0120). Zhen Wang is supported by Hunan Provincial Innovation Foundation for Postgraduate (CX20230616).
This work was carried out in part using computing resources at the High Performance Computing Platform of Xiangtan University.

\bibliographystyle{elsarticle-num} 
\bibliography{ref}

\end{document}